\documentclass[a4paper,11pt,openany]{article}
\usepackage{indentfirst}
\usepackage[iso-8859-7]{inputenc}
\usepackage{amsmath,amsfonts,amssymb}
\usepackage{graphicx}

\usepackage{color,listings,titlesec}
\usepackage[headings]{fullpage}
\usepackage{amsfonts}
\usepackage{fancyhdr}
\usepackage{graphicx}
\usepackage{amsmath}
\usepackage{marvosym}
\usepackage{hyperref}
\usepackage{amsthm}
\newcommand{\g}{\mathfrak{g}}

\newcommand{\m}{\mathfrak{m}}

\newcommand{\al}{\alpha}
\newcommand{\be}{\beta}
\newcommand{\ga}{\gamma}
\newcommand{\ad}{\mathrm{ad}}
\newcommand{\ev}{\mathrm{ev}}
\newcommand{\p}{\mathfrak{p}}
\newcommand{\di}{\mathrm{d}}
\newcommand{\su}{\mathfrak{q}}
\newcommand{\om}{\overline{\omega}}

 \def \g{\mathfrak{g}}
 
 \def \m{\mathfrak{m}}

 \def \p{\mathfrak{p}}
 \def \q{\mathfrak{q}}
 
 \def \k{\mathfrak{k}}

 \def \ad {\mathop {\hbox{\rm ad}}}
 \def \span {\mathop {\hbox{\rm span}}}

\title{$W-$ algebras and Duflo Isomorphism.}
\author{P. Batakidis\footnote{Department of Mathematics, Penn State University, USA. E-mail: batakidis@psu.edu}\;\; \& \;N. Papalexiou\footnote{Department of Mathematics, University of the Aegean, Greece. E-mail: papalexi@aegean.gr}}

\begin{document}
\maketitle
\begin{abstract}
\noindent We prove that when Kontsevich's deformation quantization is applied on weight homogeneous Poisson structures,  the operators in the $\ast-$ product formula are weight homogeneous. In the linear Poisson case $X=\g^\ast$ for a semi simple Lie algebra $\g$. As an application we provide an isomorphism between the Cattaneo-Felder-Torossian reduction algebra $H^0(\g,\m,\chi)$ and the $W-$ algebra $(U(\g)/U(\g)\m_\chi)^\m$. We also show that in the $W-$ algebra setting, $(S(\g)/S(\g)\m_\chi)^\m$ is polynomial. Finally, we compute generators of $H^0(\g,\m,\chi)$ as a deformation of $(S(\g)/S(\g)\m_\chi)^\m$.
\vspace{0.5cm}

\end{abstract}

\section{Introduction}
\noindent\textbf{1.1} In \cite{K} Kontsevich solved the deformation quantization problem of Poisson manifolds, proving the Formality Theorem for the $L_{\infty}-$ algebras $\mathcal{T}_{poly}(\mathbb{R}^k)$ of polyvector fields,  and  $\mathcal{D}_{poly}(\mathbb{R}^k)$ of polydifferential operators of bounded order, on $\mathbb{R}^k$. The theorem states that the map \, $\mathcal{U}:\;\mathcal{T}_{poly}(\mathbb{R}^k)\longrightarrow \mathcal{D}_{poly}(\mathbb{R}^k)$ defined by its Taylor coefficients
\begin{equation}\label{formality}
\mathcal{U}_n:=\sum_{\overline{m}\geq 0}\left(\sum_{\Gamma\in\mathbf{Q_{n,\overline{m}}}}\omega_{\Gamma}B_{\Gamma}\right)
\end{equation}
is an $L_{\infty}-$ morphism and a quasi-isomorphism. Properties of this map prove that there is a bijection between the gauge equivalence classes of $\ast-$ products on $C^\infty(\mathbb{R}^k)$ and the gauge equivalence classes of Poisson structures $\pi$ on $\mathbb{R}^k$. As a consequence, (\ref{formality}) provides an explicit formula of the $\ast-$ product, denoted by $\ast_K$, associated to a Poisson structure. In particular, choosing a Poisson structure $\pi$ on $\mathbb{R}^k$, the operator $\ast_K: C^\infty(\mathbb{R}^k)[[\epsilon]]\times C^\infty(\mathbb{R}^k)[[\epsilon]]\longrightarrow C^\infty(\mathbb{R}^k)[[\epsilon]]$ defined for $F,G\in C^{\infty}(\mathbb{R}^k)$ by the formula
\begin{equation}\label{kontsevich star product}
F\ast_KG:=F\cdot G+\sum_{n=1}^{\infty}\epsilon^n\left(\frac{1}{n!}\sum_{\Gamma\in\mathbf{Q_{n,2}}}\omega_{\Gamma}B_{\Gamma}^{\pi}(F,G)\right)\;,
\end{equation}
is an associative product. Formula (\ref{kontsevich star product}) comes from (\ref{formality}) and the ingredients are smaller cases of the ones therein: $\epsilon$ is a deformation parameter, $\mathbf{Q_{n,2}}$ is a family of graphs, and $\omega_\Gamma$ is a real coefficient associated to a graph $\Gamma\in\mathbf{Q_{n,2}}$. This coefficient is computed as the integral of a differential form on a compactified concentration manifold $\mathcal{C}\subset (\mathcal{H}^+)^n\times\mathbb{R}^2$, where $\mathcal{H}^+$ is the hyperbolic half-plane. Finally $B^\pi_\Gamma$ is a linear bidifferential operator on $C^\infty(\mathbb{R}^k)\times C^\infty(\mathbb{R}^k)$. Details can be found in \cite{K}, \cite{CKTB}.

\noindent Cattaneo and Felder considered in \cite{CF2} the case of a coisotropic submanifold $C=\mathbb{R}^r$ of a Poisson manifold $X=\mathbb{R}^k$,  and generalized the results of \cite{K} in \cite{CF3}.  Let $\mathcal{T}(X)$ be the DGLA of multivector fields on $X$,  $I(C)\subset C^\infty(X)$ the ideal of functions vanishing on $C$ and $\mathcal{A}=\Gamma(C,\wedge T_X)$. The Relative Formality Theorem states that there is an $L_{\infty}-$ quasi-isomorphism from $\mathcal{T}(X,C):=lim_{\leftarrow} \mathcal{T}(X)/I(C)^n\mathcal{T}(X)$, the DGLA of multivector fields on an infinitesimal neighbourghood of C, to $\tilde{\mathcal{D}}(\mathcal{A})=\oplus_n\tilde{\mathcal{D}}^n(\mathcal{A})$
where $\tilde{\mathcal{D}}^n(\mathcal{A}):=\prod_{p+q-1=n}\mathrm{Hom}^p(\otimes^q\mathcal{A},\mathcal{A})$. In this deformation problem, the associativity is controlled by a curved $A_{\infty}-$ algebra. This $A_\infty-$ algebra is flat in the linear Poisson case $X=\g^\ast$, the dual of a Lie algebra $\g$. Its $0-$th cohomology then accepts an associative product, called the \textsl{Cattaneo-Felder} product $\ast_{CF}:\;C^{\infty}(\mathbb{R}^r)[[\epsilon]]\times C^{\infty}(\mathbb{R}^r)[[\epsilon]]\longrightarrow C^{\infty}(\mathbb{R}^r)[[\epsilon]]$. For $F,G\in C^\infty(\mathbb{R}^r)$, the $\ast_{CF}-$ product is given by the formula 
\begin{equation}\label{cf}
F\ast_{CF} G:=F\cdot G+\sum_{n=1}^{\infty}\epsilon^n\left(\frac{1}{n!}\sum_{\Gamma\in\mathbf{Q^{(2)}_{n,2}}}\omega_{\Gamma}B_{\Gamma}^{\pi}(F,G)\right).
\end{equation}
 In (\ref{cf}), $\mathbf{Q^{(2)}_{n,2}}$ is a family of graphs with two \textit{colors}, a notion to be explained in $\mathcal{x}$ 2.1. The coefficient $\omega_{\Gamma}$ and the bidifferential operator $B_{\Gamma}^\pi$ are  computed similarly to the case $C=X$. 

\noindent The results of \cite{CF2} where applied by Cattaneo and Torossian in \cite{CT} to the case of symmetric spaces. In particular, let $\m\subset\g$ be a Lie subalgebra, $\chi$ a character of $\m$ and $\su$ a supplementary of $\m$ in $\g$ forming a symmetric pair. The authors considered the biquantization problem for $X=\g^\ast$ and $C=\chi+\m^\bot$. Among the basic objects of study were the aforementioned $0-$th cohomology $H^0(\m^\bot,\chi,\su,\epsilon)$, called $\epsilon-$ reduction algebra, and the reduction algebra $H^0(\m^\bot,\chi,\su)$, associated to the data $\g,\m,\su,\chi$. These two cohomology spaces are different, see Section 1.3 below. Based on that work, the main result of \cite{BAT3} is that $H^0(\m^\bot,\chi,\su,\epsilon)$ is isomorphic to $\left(U_{(\epsilon)}(\g)/U_{(\epsilon)}(\g)\m_{\chi+\rho}\right)^{\m}$ for any Lie algebra $\g$, subalgebra $\m$ and character $\chi$. Here $\rho(H)=\omega_{\Gamma'}\mathrm{Tr}(\mathrm{ad}H)$, where $\Gamma'$ is a short loop (Figure 3.1 in \cite{BAT}) with $\omega_{\Gamma'}=\frac{1}{2}$.\newline
\noindent\textbf{1.2} The systematic study of $W-$ algebras began with the paper of Premet \cite{PRE1}. The motivation is to study the finite dimensional irreducible representations of the universal enveloping algebra $U(\g)$ of a semisimple Lie algebra $\g$. Using the 1-1 correspondence of such representations with the primitive ideals of $U(\g)$, the standard approach is to study the  finite dimensional irreducible representations of the $W-$ algebra and then pass the results on $U(\g)$ via Skryabin's equivalence (Appendix in \cite{PRE1}). We first review the construction of the $W-$ algebra used by Premet (see \cite{PRE1}-\cite{PRE2}). Fix a nilpotent element $e\in\g$ and pick $h,f\in\g$ forming an $\mathfrak{sl}_2$- triple with $e$. There exists a $\g-$ invariant bilinear form $(\cdot,\cdot)$ such that $(e,f)=1$. Let $\chi\in\g^\ast$ be defined by $\chi(x)=(e,x)$ for all $x\in\g$. Set  $\g(i):=\{\xi\in\g|[h,\xi]=i\xi\}$ to be the eigenspaces of the $\mathrm{ad}h-$ action. Consider the skew-symmetric form $\omega_{\chi}$ on $\g$ defined by $\omega_{\chi}(\xi,\eta)=\chi([\xi,\eta])$. The restriction $\omega_{\chi}|_{\g(-1)}$ is non-degenerate so one can pick a lagrangian subspace $\mathfrak{l}\subset\g(-1)$. Set $\mathfrak{m}:=\mathfrak{l}\bigoplus\oplus_{i\leq-2}\g(i)$ (so $\chi$ is a character of $\mathfrak{m}$). Set $\mathfrak{m}_{\chi}$ to be the space generated by the elements $\{\xi-\chi(\xi),\;\xi\in\mathfrak{m}\}$ and let $\su$ be such that $\g=\m\oplus\su$. The space $\mathcal{S}:=e+\mathrm{ker}\ad(f)$ is then called the \textsl{Slodowy slice} at $e$ through the adjoint orbit $\mathcal{O}:=G\cdot e$, (see \cite{PRE1} and \cite{GaGi} $\mathcal{x}$ 1.2). The (finite) \textsl{$W$- algebra} corresponding to the data $(\g,e)$ is then defined as $U(\g,e)=\text{End}_{\g}(Q_{\chi})^{\text{op}}$, where $Q_{\chi}$ is the $U(\g)$-module induced from the one dimensional left $U(\m)$-module obtained by the character $\chi$. From the PBW theorem follows that $U(\su)\simeq U(\g)/U(\g)\mathfrak{m}_{\chi}$ and $U(\g,e)$ can then be defined equivalently  as the quantum Hamiltonian reduction $(U(\g)/U(\g)\mathfrak{m}_{\chi})^{\mathfrak{m}}$. The associated graded algebra $\mathrm{\textbf{gr}}U(\g,e)$ is isomorphic to the graded algebra of functions on  $\mathcal{S}$. 
At the same time, Losev's approach in \cite{L2},\cite{L1} is close to deformation quantization and begins with a more geometric definition of the $W$- algebra as the, specialized at $\hbar=1$, $G-$ invariants of the $\mathbb{K}[\hbar]-$ algebra $\mathbb{K}[X][\hbar]$, where $X=G\times\mathcal{S}$. This algebra is equipped with Fedosov's $\ast-$ product (\cite{FED}).

\noindent\textbf{1.3} We use the approach of Kontsevich and Cattaneo-Felder-Torossian on deformation quantization to prove that when $\mathbb{R}^k$ carries a weight homogeneous Poisson structure with respect to a weight vector $\om$, the terms of the $\ast-$ product (\ref{kontsevich star product}) applied to weight homogeneous functions $F,G\in C^\infty(\mathbb{R}^k)$, are weight homogeneous with respect to $\om$. As a consequence, in the semisimple Lie algebra case, the degree of the second term in the $\ast-$ product $\ast_{CF}$ recovers the quasihomogeneous degree of the transverse Poisson structure on a slice, e.g $\mathcal{S}$, see \cite{DSV}. More precisely, if $x_i, \;i=1,\ldots,k$ is a basis of $\g$ and $n_i$ are the weights of the $\mathrm{ad}h-$ action, $[h,x_i]=n_ix_i$, set $\om_\g=(n_1+2,\ldots,n_k+2)$. The weight vector $\om_{\g}$ is the one inducing the Kazhdan grading on the symmetric algebra $S(\g)$ (Remark \ref{kazhdan weight cor}). The fact that the quasihomogenous degree of the Lie-Poisson structure in this case is $-2$, was proved in \cite{DSV} $\mathcal{x} 3$, while the Slodowy slice case was also implicit in \cite{GaGi}. Our first main result, (see Theorem \ref{kazhdan theorem} for the proof), states that when $\g$ is semisimple and the weight is $\om_\g$, there is an isomorphism between the reduction algebra $H^0(\m^{\bot},\chi,\su)$ and $U(\g,e)$, providing a new model of the $W-$ algebra:
\newtheorem{con}{Theorem}[section]
\begin{con}\label{fon}
Let $\g$ be a semisimple Lie algebra, and $\{e,h,f\}$ an $\mathfrak{sl}_2-$ triple. Let $\chi,\m_{\chi}$ be defined from the $W-$ algebra construction. There is an associative algebra isomorphism
\begin{equation}
\mathcal{Q}:\; H^0(\m^\bot,\chi,\su)\stackrel{\simeq}{\longrightarrow} \left(U(\g)/U(\g)\m_{\chi}\right)^{\m}
\end{equation}
\end{con}

 This way, we transfer the study of $W$-algebras to the study of $H^0(\m^\bot,\chi,\su)$. When $\mathcal{O}=G\cdot e$ is principal, Theorem \ref{fon} recovers the Duflo algebra isomorphism $S(\g)^G\stackrel{\sim}{\longrightarrow} \mathcal{Z}(\g)$ between the $G-$ invariants of $S(\g)$ and the center of $U(\g)$, as in \cite{K} pp. 207-213. Note that $H^0(\m^{\bot},\chi,\su)$ is different from  $H^0(\m^\bot,\chi,\su,\epsilon)$. In fact, for every Lie algebra $\g$, subalgebra $\m\subset\g$ and character $\chi$ of $\m$, it is $H^0(\m^\bot,\chi,\su,\epsilon=1)\hookrightarrow H^0(\m^\bot,\chi,\su)\hookrightarrow\left(U(\g)/U(\g)\m_{\chi+\rho}\right)^{\m}$ (see Lemma 3.1 and Proposition 3.4 of \cite{BAT}). In the $W$-algebra setting the character $\rho\in\m^\ast$ is missing, since $\m$ is a nilpotent subalgebra and so $\rho(H)=0,\;\forall H\in\m$ by Engel's Theorem. Hence we have the direction $H^0(\m^\bot,\chi,\su)\hookrightarrow\left(U(\g)/U(\g)\m_{\chi}\right)^{\m}$. We prove the inverse direction  using the grading aspect of the system of equations (\ref{reduction equations 3}) defining $H^0(\m^{\bot},\chi,\su)$. As we show in Section 3, these are now homogeneous, and provide a certain filtration in $H^0(\m^{\bot},\chi,\su)$ (Remark \ref{degree in H}). In Section 4 we  prove that, if $\g_e$ is the centralizer of $e$ in $\g$, each element of $H^0(\m^{\bot},\chi,\su)$ is uniquely determined by an element of $S(\g_e)$. Finally we compute precisely the generators of $H^0(\m^{\bot},\chi,\su)$  (see Theorem \ref{paragetai} for the proof); if $x_1,\ldots,x_r$ is a basis of $\g_e$, we construct elements $\tilde{P}_1,\ldots, \tilde{P}_r\in H^0(\m^{\bot},\chi,\su)$ such that 
\newtheorem{cre}[con]{Theorem}
\begin{cre}\label{t}
The algebra $ H^0(\m^\bot,\chi,\su)$ is generated by the elements $\tilde{P}_i,\;i=1,\ldots,r$.
\end{cre}

Our results are closely related with the polynomial conjecture proposed by C. Torossian in \cite{Toro}. More precisely, if $(\g,\sigma)$ is a symmetric pair and $\g=\k\oplus\p$ is the decomposition relative to $\sigma$, the conjecture suggests that $(U(\g)/U(\g)\k)^{\k}$ and $S(\p)^{\k}$ are isomorphic as algebras. If $\g$ is semisimple, and $\m,\chi$ are given by the $W$-algebra construction, we prove that there is an algebra isomorphism between $(U(\g)/U(\g)\m_{\chi})^{\m}$ and a deformation of $(S(\g)/S(\g)\m_{\chi})^{\m}$.
\newpage

\section{Deformation Quantization and Weight homogeneous Poisson structures.}
\noindent\textbf{2.1. Deformation Quantization background.}\newline

\noindent \textbf{2.1.1. Some notation.} Let $\mathbb{K}$ be a field of characteristic zero containing $\mathbb{R}$, and $\g$ a Lie algebra of finite dimension over $\mathbb{K}$. Let $\m\subset\g$ be a subalgebra, $\chi$ a character of $\m$,  $S(\g)$ the symmetric and $U(\g)$ the universal enveloping algebra of $\g$, respectively. Set $\m_{\chi}$ to be the vector subspace of $S(\g)$ generated by the set $\{m-\chi(m),\;m\in\m\}$ and denote as $S(\g)\m_\chi$, $U(\g)\m_{\chi}$ the ideal of $S(\g)$ and right ideal of $U(\g)$, respectively, generated by $\m_{\chi}$. There is a natural isomorphism between $S(\g)$ and $\mathbb{K}[\g^\ast]$, the algebra of polynomials on the dual Lie algebra $\g^\ast$. The algebra $S(\g)$ is equipped with a natural Poisson structure defined for $x_1,x_2\in\g$ by $\{x_1,x_2\}:=[x_1,x_2]$ turning $\g^\ast$ into a Poisson manifold. Furthermore, the algebra $(S(\g)/S(\g)\m_{\chi})^{\m}$, of $\ad\m-$ invariants inherits a Poisson structure. Let then $\m^{\bot}:=\{l\in\g^{\ast}/l(\m)=0\}$ and $\mathbb{K}[\chi+\m^{\bot}]^{\m}$ be the Poisson algebra of $\m$-invariant polynomial functions on $\chi +\m^{\bot}$. One then has $(S(\g)/S(\g)\m_{\chi})^{\m}\simeq \mathbb{K}[\chi+\m^{\bot}]^{\m}$ as algebras. Considering $\m_{\chi}^{\bot}:=\{f\in\g^{\ast}/f|_{\m}=\chi\}$ as a coisotropic submanifold of $\g^\ast$ it is possible to apply the biquantization techniques of \cite{CT} to write the corresponding $\ast_{CF}-$ product for this algebra; we briefly recall some of the necessary definitions adjusted to our setting.\newline

\noindent \textbf{2.1.2. Kontsevich's construction.}  
\noindent Denote by $\mathbf{Q_{n,2}},\;n\in\mathbb{N}$ the set of all \textsl{admissible} graphs $\Gamma$, meaning graphs with the following properties: The set $V(\Gamma)$ of vertices of $\Gamma$ is the disjoint union of two ordered sets $V_1(\Gamma)$ and $V_2(\Gamma)$, isomorphic to $\{1,\ldots,n\}$ and $\{1,2\}$ respectively. Their elements are called \textsl{type I} vertices, for $V_1(\Gamma)$, and \textsl{type II} vertices, for $V_2(\Gamma)$.  The set $E(\Gamma)$ of edges in the graph is finite. Each edge starts from a type I vertex and ends to a vertex in $V_1(\Gamma)\cup V_2(\Gamma)$ (no loops or double edges). A vertex $\upsilon\in V_1(\Gamma)$ receiving no edge will be called a \textsl{root}. All elements of $E(\Gamma)$ are oriented and the set of edges $S(r)$ starting from $r\in V_1(\Gamma)$ is ordered. This induces an order on $E(\Gamma)$, compatible with the order on $V_1(\Gamma)$ and $S(r)$.

 \noindent Consider a Poisson structure $\pi$ on $\mathbb{R}^k$. To a graph $\Gamma\in \mathbf{Q_{n,2}}$ one associates a bidifferential operator $B_{\Gamma}$ as follows: Let $\{x_1,\ldots,x_k\}$ coordinate functions on $\mathbb{R}^k$ and $L:\;E(\Gamma)\longrightarrow\{x_1,\ldots,x_k\}$ be a labeling function for the edges of $\Gamma$. Fix a vertex $r\in\{1,\ldots,n\}$. If $\mathrm{card}(S(r))\neq 2$, set $B_{\Gamma}=0$\footnote{This is actually an implication of the general construction.}. If $\mathrm{card}(S(r))= 2$, let $S(r)=\{e_r^1,e_r^2\}$ be the ordered set of edges leaving $r$. Associate the bracket $\{L(e_r^1),L(e_r^2)\}$ to $r$. To each vertex $1,2\in V_2(\Gamma)$ associate respectively a function $F,G\in S(\g)$, and to the $p^{th}-$ edge of $S(r)$, associate the partial derivative with respect to the coordinate variable $L(e_r^p)$. This derivative acts on the function associated to $v\in V(\Gamma)$ where the edge $e_r^p$ arrives.
Since $E(\Gamma)\subset V_1(\Gamma)\times (V_1(\Gamma)\cup V_2(\Gamma))$, let $(p,m)\in E(\Gamma)$ represent an oriented edge of $\Gamma$ from $p$ to $m$. Then define
\begin{equation}\label{bidiff op}
B_{\Gamma}^\pi(F,G)= \sum_{L:E(\Gamma)\rightarrow \{1,\ldots,k\}}\left[\prod_{r=1}^{\#(V_1(\Gamma))}\left(\prod_{\delta\in E(\Gamma),\;\delta=(\cdot,r)}\partial_{L(\delta)}\right)\{L(e_r^1),L(e_r^2)\}\right]\times
\end{equation}
\[\times\left(\prod_{\delta\in E(\Gamma),\;\delta=(\cdot,1)}\partial_{L(\delta)}\right)(F)\times\left(\prod_{\delta\in E(\Gamma),\;\delta=(\cdot,2)}\partial_{L(\delta)}\right)(G).\] We drop the exact definition of the coefficient $\omega_{\Gamma}$ in (\ref{kontsevich star product}) since it can be found in the given references. In deformation quantization of a Poisson manifold, one trivially has a single choice of the \textsl{color} of variables available for every edge in a graph $\Gamma\in\mathbf{Q_{n,2}}$ since for  $e\in E(\Gamma)$, $L(e)$ determines necessarily a coordinate of $\mathbb{R}^k$. In biquantization we consider two \textsl{colors}; if $\g$ is a Lie algebra and $\m$ a subalgebra, we discuss briefly here the case $X=\g^\ast$ and $C=\m^\bot$ for later use, but the same statements hold for $X=\mathbb{R}^k,\;C=\mathbb{R}^r$. Suppose $\su$ is a supplementary of $\m$, i.e $\g=\m\oplus\su$, $\{m_1,\ldots, m_t\}$ is a basis of $\m$ and $\{q_1,\ldots, q_r\}$ a basis of $\su$. We identify spaces $\su^{\ast}\simeq\g^{\ast}/\m^{\ast}\simeq\m^{\bot}$.  For $e\in E(\Gamma)$, let $c_e\in\{+,-\}$ be its color defined as follows: Consider a 2-color label $L:\;E(\Gamma)\longrightarrow \{m_1,\ldots,m_t,q_1,\ldots,q_r\}$, satisfying $L(e)\in\{m_1,\ldots,m_t\}$ if $c_e=-$ and $L(e)\in\{q_1,\ldots,q_r\}$ if $c_e=+$. This way, the dual variables $\{m_1^{\ast},\ldots, m_t^{\ast}\}$ of $\m^{\ast}$ are associated to the color $(-)$ and the variables $\{q_1^{\ast},\ldots, q_r^{\ast}\}$ of $\su^{\ast}$ are associated to $(+)$. Graphically, the color $(-)$ will be represented with a dotted edge and the color $(+)$ will be represented with a solid edge. The corresponding formulas (\ref{bidiff op}) and  (\ref{kontsevich star product}) need some modifications in biquantization; for $F,G\in S(\su)$, one has to use the 2-colored label $L$ that we just described. From now on, all graphs, their associated operators $B_\Gamma$ and coefficients $\omega_\Gamma$ are colored. We denote by $\mathbf{Q_{n,2}^{(2)}}$ the family of admissible graphs with $2$ colors\footnote{Double edges are not allowed, meaning edges with the same color, source and target.}. \newline

\noindent\textbf{2.2. Weight homogeneous Poisson structures.}\newline

This Section concerns the weights in deformation quantization when applied to weight homogeneous Poisson structures. This calculus can be generalized to the Formality Theorem if one associates  weight homogeneous, with respect to the same weight vector, polyvector fields $P_1,\ldots,P_n$ to the type I vertices of a $\Gamma\in \mathbf{Q_{n,\overline{m}}}$ . The notation used is from \cite{LPV}, $\mathcal{x}$ 8.1.3.\newline

\noindent\textbf{2.2.1. Definitions.} Consider $\mathbb{R}^k$ with coordinates $x_1,\ldots,x_k$. Let $\overline{\omega}=(\overline{\omega}_1,\ldots,\om_k)$ be a $k-$ tuple of positive integers. An $F\in C^\infty(\mathbb{R}^k)$ is called \textsl{weight homogeneous} with respect to $\om$ if there is an $r\in\mathbb{N}$ such that $F(\lambda^{\om_1}x_1,\ldots,\lambda^{\om_k}x_k)=\lambda^rF(x_1,\ldots,x_k)$, $\forall \lambda\in \mathbb{R}$. In this case, $\om$ is  called a \textsl{weight vector} and the number $r$ is called the \textsl{weight} of $F$. We then write $\om(F)=r$. Similarly, a $p-$ vector field $P$ on $\mathbb{R}^k$ is said to be \textsl{weight homogeneous} with respect to $\om$ if applying $P$ to weight homogeneous functions $F_1,\ldots,F_p\in C^\infty(\mathbb{R}^k)$ we get a weight homogeneous smooth function. It turns out that if $P,F_1,\ldots,F_p$ are weight homogeneous of weights $\om(P),\om(F_1),\ldots,\om(F_p)$ respectively, then $P(F_1,\ldots,F_p)$ is weight homogeneous of weight $\om(P)+\sum_{i=1}^n\om(F_i)$ (or $P(F_1,\ldots,F_p)=0$). The weighted Euler vector field $E_{\om}=\sum_{i=1}^k\om_ix_i\frac{\partial}{\partial x_i}$, traces the weights of homogeneous elements; $E_{\om}(F)=\om(F)F$ if and only if $F$ is weight homogeneous. For weight homogeneous polyvector fields $P$  one gets the corresponding result applying the Lie derivative $\mathcal{L}_{E_{\om}}$ to $P$: $\mathcal{L}_{E_{\om}}P=\om(P)P$.

\noindent\textbf{2.2.2} \textbf{Weight homogeneous degree for $\ast-$ products.} Let $\pi$ be a weight homogeneous Poisson structure on $\mathbb{R}^k$. Equivalently we think of it as a weight homogeneous bivector field satisfying $[\pi,\pi]_{\mathrm{SN}}=0$ for the Schouten-Nijenhuis bracket. The next Lemma computes the weight homogeneous degrees of the terms in $\ast_{CF}$ for the case of a coisotropic submanifold $C=\mathbb{R}^r$ of $X=\mathbb{R}^k$. 

\newtheorem{cet}[con]{Lemma}
\begin{cet}\label{kazhdan degree product}
Let  $\pi$ a weight homogeneous Poisson structure on $\mathbb{R}^k$ and $F,G$ weight homogeneous functions in $C^\infty(\mathbb{R}^{r})$ with respect to a weight vector $\om$. If $\Gamma\in\mathbf{Q^{(2)}_{n,2}}$, then
\begin{equation}\label{kazhdan product}
\om(B_{\Gamma}(F,G))=
\om(F)+\om(G)+n\cdot\om(\pi)
\end{equation}
\end{cet}

\begin{proof} Fix a graph $\Gamma\in\mathbf{Q^{(2)}_{n,2}}$ and a label $L$. Suppose that $\Gamma$ has with $j$ roots and set
\[P=\prod_{r=1}^{\#(V_1(\Gamma))}\left(\prod_{\delta\in E(\Gamma),\;\delta=(\cdot,r)}\partial_{L(\delta)}\right)\{L(e_r^1),L(e_r^2)\},\;\;D_1=\prod_{\delta\in E(\Gamma),\;\delta=(\cdot,1)}\partial_{L(\delta)},\;\;D_2=\prod_{\delta\in E(\Gamma),\;\delta=(\cdot,2)}\partial_{L(\delta)}\] 
so that the operator corresponding to $\Gamma$ and $L$ is, given (\ref{bidiff op}), $B_\Gamma=PD_1D_2$.
Then
\begin{equation}\label{imi} 
\om(B_{\Gamma}(F,G))=\om(P)+\om(D_1F)+\om(D_2G).\end{equation}
We have 
\[\om(P)=\sum_{\substack{r=1\\\text{r is root}}}^j\om(\{L(e_r^1),L(e_r^2)\})+
\sum_{\substack{r=1\\ r\;\;\text{is not root}}}^{n-j}\om(\partial_{L(e_{r_1})\cdots L(e_{r_s})}\{L(e_r^1),L(e_r^2)\})\]
and
\[\om(D_1F)=\om(F)
-\sum_{p=1}^{\#\rightarrow F}\om(L(e_p)).\] 
The expression $\sum_{p=1}^{\#\rightarrow F}\om(L(e_p))$ sums the weights of labels on edges pointing to $F$, i.e $p$ runs the set of vertices carrying an edge, $e_p$, towards $F$. Furthermore, $L(e_{r_1})\cdots L(e_{r_s})$ stand for the labels of edges pointing to the vertex $r\in V_1(\Gamma)$. Summing all terms, (\ref{imi}) is 
\[\om(B_{\Gamma}(F,G))=\sum_{\substack{r=1\\\text{p is root}}}^j\om(\{L(e_r^1),L(e_r^2)\})+
\sum_{\substack{r=1\\ r\;\;\text{is not root}}}^{n-j}\om(\partial_{L(e_{r_1})\cdots L(e_{r_s})}\{L(e_r^1),L(e_r^2)\})\]
\begin{equation}\label{imia}
+\om(F)
-\sum_{p=1}^{\#\rightarrow F}\om(L(e_p))+\om(G)
-\sum_{p=1}^{\#\rightarrow G}\om(L(e_p)).
\end{equation}
It is immediate by the definition that since  $\pi$ is weight homogeneous,
\begin{equation}\label{vertex weight}
\om(\{L(e_r^1),L(e_r^2)\})=\om(\pi)+\om(L(e_r^1))+\om(L(e_r^2)).
\end{equation}  This is the contribution of a root vertex $r$ to the weight of the polynomial $B_\Gamma(F,G)$. When $r\in V_1(\Gamma)$ is not a root, the weight contributed from $r$ to the weight of  $B_\Gamma(F,G)$ is \[\om(\partial_{L(e_{r_1})\cdots L(e_{r_s})}\{L(e_r^1),L(e_r^2)\})=\om(\pi)+\om(L(e_r^1))+\om(L(e_r^2))-\left(\om(L(e_{r_1}))+\cdots \om(L(e_{r_s}))\right).\] Summing weights over all type I vertices, (\ref{imia}) gives
\[\om(B_\Gamma(F,G))=\om(F)+\om(G)+n\cdot\om(\pi).\]
 \end{proof}

\section{A new $W$-algebra model.}
\label{reduction algebra}
\noindent \textbf{3.1.} \textbf{Reduction algebras.} We describe some particular types of graphs (see ~\cite{CT} $\mathcal{x}$ 1.3, and ~\cite{BAT} $\mathcal{x}$ 2.3). They are colored graphs with only one type II vertex, see Figure 1. Denote as $e_{\infty}$ an edge colored by $(-)$ with no end.  
\newtheorem{conn}[con]{Definition}
\begin{conn}\label{conn}
\begin{enumerate}
\item \textbf{Bernoulli.} Graphs of this type, with $i$ type I vertices, $i\in\mathbb{N}_{\geq 2}$, will be denoted by $\mathcal{B}^i$. They have $2i$ edges, and $i$ of them are pointing to the type II vertex. They have an $e_\infty$ edge and a root.
\item \textbf{Wheels.} Graphs of this type, with $i$ type I vertices, $i\in\mathbb{N}_{\geq 2}$, will be denoted by $\mathcal{W}^i$. They have $2i$ edges and leave no edge to $\infty$. Furthemore, $i$ of the edges form an oriented polygon and the rest point to the type II vertex. 
\item \textbf{Bernoulli attached to a wheel.} Graphs of this type, with $i$ type I vertices, $i\in\mathbb{N}_{\geq 3}$, will be denoted by $\mathcal{BW}^i$. They have $i-1$ edges towards the type II vertex and leave an edge to $\infty$.
For an $\mathcal{W}^m-$ type graph $W_m$ attached to a $\mathcal{B}^l-$ type graph $B_l$, we will write $B_lW_m\in\mathcal{B}^l\mathcal{W}^m$. Obviously $\mathcal{B}^l\mathcal{W}^m\subset\mathcal{BW}^{l+m}$.
\end{enumerate}
\end{conn}
\begin{figure}[h!]
\begin{center}\label{tena}
\includegraphics[width=9cm]{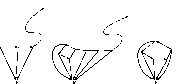}
\caption\footnotesize{From left to right, a $\mathcal{B}^3$-type graph, a $\mathcal{B}^3\mathcal{W}^4$-type graph, and a $\mathcal{W}^5$-type graph.}
\end{center}
\end{figure}

Set $\mathbf{Q^{\infty}_{i,1}}$ to be the family of such graphs with $i$ type I vertices, namely graphs of the categories 1 and 3 of the previous  definition. Consider $X=\g^\ast$, the dual of a Lie algebra $\g$ endowed with the Poisson structure $\{\cdot,\cdot\}$. Let $\m\subset\g$ be a subalgebra, $\chi$ a character of $\m$ and $\su$ such that $\g=\m\oplus\su$. Fix a basis $\{m_1,\ldots,m_t\}$ of $\m$. Let $\{e^1_l,e^2_l\}$ be the ordered set of edges leaving $l\in V_1(\Gamma)$ of a colored graph $\Gamma \in\mathbf{Q^{\infty}_{i,1}}$. Let $B_{\Gamma}:\;S(\su)\longrightarrow S(\su)\otimes \m^{\ast}$ be the differential operator: $F\mapsto (m_j\mapsto B_{\Gamma}(m_j)(F))$ for   $j=1,2,\ldots,t$, where 
 \begin{equation}\label{form}
 B_{\Gamma}(m_j)(F)= \sum_{\substack{\text{L: colored}\\\mathrm{L}(e_{\infty})=m_j\\}}\left[\prod_{r=1}^{i}\left(\prod_{e\in E(\Gamma),\;e=(\cdot,r)}\partial_{L(e)}\right)\{L(e_r^1),L(e_r^2)\}\right]\times \left( \prod_{\substack{e\in E(\Gamma)\\ e=(\cdot,F)\\}}\partial_{L(e)}F\right) 
\end{equation}
 For $m_j\in\{m_1,\ldots,m_t\}$, we define $\di^i:\;S(\su)[\epsilon]\longrightarrow S(\su)[\epsilon]\otimes \m^{\ast}$ to be the differential operator
\[F\mapsto (m_j\mapsto \di^i(m_j)(F)),\]
where \[\di^i(m_j)(F)=\sum_{\Gamma\in \mathbf{Q^\infty_{i,1}}}\omega_{\Gamma}B_{\Gamma}(m_j)(F).\]
Set finally $\di(\m^\bot,\chi,\su,\epsilon):\;S(\su)[\epsilon]\longrightarrow S(\su)[\epsilon]\otimes \m^{\ast}$ be the differential operator

 \begin{equation}\label{di}
\di(\m^\bot,\chi,\su,\epsilon)=\sum_{i=1}^{\infty}\epsilon^i \di^i
\end{equation}

One may repeat the same construction without $\epsilon$ to define an operator $\di(\m^\bot,\chi,\su):\;S(\su)\longrightarrow S(\su)\otimes\m^{\ast}$ with $\di(\m^\bot,\chi,\su)=\sum_{i=1}^{\infty} \di^i$.

\newtheorem{new}[con]{Remark}
\begin{new}\label{inv}
A direct computation shows that $\di^1(F)=0$ if and only if $F\in S(\su)^\m$. 
\end{new}

\newtheorem{demek}[con]{Definition}
\begin{demek}(\cite{CT})\label{redalg}
Consider $X=\g^\ast$ with the (Lie-) Poisson structure and $C=\chi+\m^\bot$ as a coisotropic submanifold.\newline
a) The $\mathbf{\epsilon}$-\textbf{reduction algebra} over $\chi+\m^\bot$ is the vector space of solutions $F_{(\epsilon)}\in S(\su)[\epsilon]$ of the equation

\begin{equation}\label{reduction equations 2}
\di(\m^\bot,\chi,\su,\epsilon)(F_{(\epsilon)})=0
\end{equation}

 equipped with the $\ast_{CF}-$ product (\ref{cf}). We denote this algebra as $H^0(\m^\bot,\chi,\su,\epsilon)$. \newline

\noindent b)  The \textbf{reduction algebra} over $\chi+\m^\bot$ is the vector space of solutions $F\in S(\su)$~of the equation

\begin{equation}\label{reduction equations 3}
\di(\m^\bot,\chi,\su)(F)=0
\end{equation}

equipped with the $\ast_{CF}$-product (without $\epsilon$). We denote this algebra as  $H^0(\m^\bot,\chi,\su)$.
\end{demek}
For simplicity we denote as $\di=\sum_{i=1}^\infty \di^i$ the differential in both the defining equations (\ref{reduction equations 2}) and (\ref{reduction equations 3}) since the existence or not of the deformation parameter $\epsilon$ will be explicitly indicated when needed. In fact, when $\g$ is semisimple,  we will show in Proposition \ref{homo} that the two systems of differential equations (\ref{reduction equations 2}) and (\ref{reduction equations 3})  are in a certain sense equivalent. 

\noindent Let $\mathrm{deg}_{\su}(G)=p$ mean that $G\in S(\su)$ has polynomial degree $p$ . Similarly, $\mathrm{deg}_{\epsilon}$ is the $\epsilon-$ degree of elements in $S(\su)[\epsilon]$ and the differential operator $\di$ on $S(\su)[\epsilon]$ as defined in (\ref{di}) . For $F\in S(\su)[\epsilon]$, set $\mathrm{deg}_{\su,\epsilon}(F):=\mathrm{deg}_{\su}(F)+\mathrm{deg}_{\epsilon}(F)$.  Let now $F_{(\epsilon)}=F_0+\epsilon F_{1}+\cdots +\epsilon^n F_n\in H^0(\m^\bot,\chi,\su,\epsilon)$.
Grouping the terms of the left-hand side of (\ref{reduction equations 2}) with respect to their $\deg_\epsilon$ we get the following  system of $\deg_\epsilon-$ homogeneous linear partial differential equations,
\begin{equation}\label{system}
\forall p\in\mathbb{N},\;\;\sum_{i=1}^{p}\di^i(F_{p-i})=0.
\end{equation}
By \cite{CT}, $\mathcal{x}$ 2, Lemma 7, only the colored graphs $\Gamma\in\mathbf{Q^\infty_{2n+1,1}}$, have a non-zero contribution to $\di$. So (\ref{system}) is
\begin{equation}\label{system2}
\forall p\in\mathbb{N}_0,\;\;\sum_{i=0}^p\di^{2i+1}(F_{2(p-i)})=0.
\end{equation}
This means that every element of $H^0(\m^\bot,\chi,\su,\epsilon)$ can be written as $F_{(\epsilon)}=\sum_{i=0}^n\epsilon^{2i}F_{2i}$. Turning equation (\ref{reduction equations 2}) into a homogeneous system is possible using the $\epsilon-$ degree and thus works for well in this case. System (\ref{reduction equations 3}) is much more complicated.

\noindent \textbf{3.2. Weights in the reduction algebra.} Consider $X=\g^\ast$ endowed with a weight homogeneous Poisson structure. Then the algebra of smooth functions on the coisotropic submanifold $\m_\chi^\bot$ is isomorphic to $S(\su)$. For a weight homogeneous $F\in S(\su)$ and a graph $\Gamma\in\mathbf{Q^\infty_{n,1}}$, we want to prove that $B_{\Gamma}(m_j)(F)$ is a weight homogeneous function and then compute its weight. We say that $B_{\Gamma}(m_j)$ is weight homogeneous, if for any weight homogeneous function $F$, the function $B_{\Gamma}(m_j)(F)$ is weight homogeneous. In that case we set $\om(B_{\Gamma}(m_j)):=\om(B_{\Gamma}(m_j)(F))-\om(F)$ to be the weight of $B_{\Gamma}(m_j)$. For the weight $\om(\pi)$ of the bivector field $\pi$ see the relative discussion at Section 2.2.1 for polyvector fields. The following Lemma computes $\om(B_{\Gamma}(m_j))$ for $\Gamma\in\mathbf{Q^\infty_{n,1}}$.

\newtheorem{conana}[con]{Lemma}
\begin{conana}\label{weight diff lemm}
 Let $\pi$ be a weight homogeneous Poisson structure on $\g^\ast$. Let $\Gamma\in\mathbf{Q^\infty_{n,1}}$, $L$ be a label function and fix $L(e_\infty)=m_j$.   Then $B_{\Gamma}(m_j)$ is weight homogeneous and
\begin{equation}\label{weight differential formula}
\om(B_{\Gamma}(m_j))=n\cdot\om(\pi)+\om(m_j)
\end{equation}
\end{conana}
\begin{proof}  Let $\Gamma\in \mathcal{B}^n$. Label as $1$ the root of the graph, as $2$ the type I vertex receiving the edge starting from $1$ and etc. Thus the $n^\text{th}$ vertex is the origin of $e_\infty$. Let $e_i^1,\;i=1,\ldots,n$ be the edges deriving $F$,  $e_i^2,\;i=1,\ldots,n-1$ be the edges between the type I vertices, and $e_n^2=e_{\infty}$.  Fixing a label $L$,  one has \[\om(\{L(e_i^1),L(e_i^2)\})=\om(L(e_{i-1}^2))\] for $i=2,\ldots,n$. 
Borrowing arguments from Lemma \ref{kazhdan degree product}, 
\[\om(B_{\Gamma}(m_j)(F))=\om(F)-\sum_{i=1}^n\om(L(e_i^1))+
\om(\{L(e_1^1),L(e_1^2)\}).\]
Analyzing the weight $\om(\{L(e_1^1),L(e_1^2)\})$,  we get 
\[\om(\{L(e_1^1),L(e_1^2)\})=\om(F)-\sum_{i=1}^n\om(L(e_i^1))+\om(\{L(e_1^1),\ldots,\{L(e_{n-2}^1),\{L(e_{n-1}^1),\{L(e_n^1),L(e_n^2)\}\cdots\}).\]
Finally, by applying relation (\ref{vertex weight}) to the right hand side of the last equation, we conclude that
\[\om(B_{\Gamma}(m_j)(F))=\om(F)+\om(m_j)+n\cdot\om(\pi).\]

\noindent The case $\Gamma\in\mathcal{BW}^n$ works similarly: Let $p$ type I vertices be in the wheel and $n-p$ type I vertices be in the Bernoulli part of $\Gamma$. Label as $1$ the vertex inside the wheel that receives the edge leaving the vertex of the wheel where the Bernoulli part is attached. Label the rest of the vertices following the orientation of the edges in the wheel and then in the Bernoulli part such that $e_\infty$ leaves the vertex labeled as $n$.  Order the edges deriving $F$ by $e_i^1$ and  the edges among the type I vertices by $e_i^2$, $i\in\{1,\ldots,n\}$. Again $e_n^2=e_\infty$. In particular, none of the edges leaving the $p^{\text{th}}$ vertex  derives $F$, so we will name as $e_p^2$ the edge inside the wheel and as $e_p^1$ the edge towards the Bernoulli part of $\Gamma$. For $i>1$, $i\notin\{p,p+1\}$, the contribution of the $i^{\text{th}}$ type I vertex to the weight of the polynomial function $B_{\Gamma}(m_j)(F)$ is $\om(\pi)+\om(L(e_i^1))+\om(L(e_i^2))-\om(L(e^2_{i-1}))$. Similarly, at the first type I vertex, the contribution is $\om(\pi)+\om(L(e_1^1))+\om(L(e_1^2))-\om(L(e^2_p))$. Using the notation of Lemma \ref{kazhdan degree product} and formula (\ref{form}), let $B_{\Gamma}(m_j)=PD_1$. Since $\Gamma\in\mathcal{BW}^n$, $B_{\Gamma}(m_j)$ is a constant coefficient operator, that is $P$ is a constant and $\om(P)=0$. Then
\[\om(B_{\Gamma}(m_j)(F))=\om(F)-\sum_{i=1}^n\om(L(e_i^1)).\]
A similar calculation as before proves the claim. \end{proof}

\newtheorem{coni}[con]{Lemma}
\begin{coni}\label{kazhdan wheel lemma}
Fix a $\Gamma\in\mathcal{W}^n$ and let $F\in S(\su)$.  Then
\begin{equation}\label{kazhdan wheel}
\om(B_{\Gamma}(F))=\om(F)+n\cdot\om(\pi)
\end{equation}
\end{coni}
\begin{proof} One can get the claim with similar arguments as in the previous Lemma. \footnote{Otherwise the Lemma is straightforward since the wheels are graphs that fall into the range of Lemma \ref{kazhdan degree product}.} \end{proof}

The following Lemma shortens the computations in the case of a wheel-type graph.
\newtheorem{rato}[con]{Lemma}
\begin{rato}\label{k}
Let $\Gamma\in\mathcal{W}^n$ and fix a weight $\om$. Let $e_i^1$ be the edges deriving the function $F$  associated to the type II vertex. Then
\[\sum_{i=1}^n\om(L(e_i^1))=-n\cdot\om(\pi)\]
\end{rato}
\begin{proof}
Pick a type I vertex and label it as $1$. Continue labeling the type I vertices following the edges inside the wheel. At each type I vertex $i=1,\ldots,n-1$, it is  $\om(L(e_i^2))=\om(L(e_{i+1}^1))+\om(L(e_{i+1}^2))+\om(\pi)$ and $\om(L(e_n^2))=\om(L(e_1^1))+\om(L(e_1^2))+\om(\pi)$ at the $n^{\text{th}}$ vertex. One then has to sum by parts and simplify.
\end{proof}

We say that a Poisson structure is weight homogeneous if the associated bivector $\pi$  is weight homogeneous according to $\mathcal{x}$ 2.2.1.

\newtheorem{tso}[con]{Proposition}
\begin{tso}\label{homo}
Consider a weight homogeneous (Lie-) Poisson structure $\pi$ on $\g^\ast$, $\m\subset\g$ a Lie subalgebra, $\chi$ a character of $\m$ and $\su\subset\g$ a vector subspace such that $\g=\m\oplus\su$. Then the defining system (\ref{reduction equations 3}) is equivalent to (\ref{system2}).
\end{tso}
\begin{proof} We first observe that even since the result of Lemma \ref{weight diff lemm} depends on the label $L$, this will not affect our use of it in writing (\ref{reduction equations 3}) as a system of  homogeneous equations with respect to $\om$.
Let $F=\sum_{i=0}^{k_0-1}\tilde{F}_{i}\in H^0(\m^\bot,\chi,\su)$, with $\tilde{F}_i\in S(\su)$ homogeneous of $\om(\tilde{F}_i)=k_0-i$. Then  (\ref{reduction equations 3}) is 
\begin{equation}\label{prohomo}
\sum_{i=0}^{k_0-1}\di(\tilde{F}_i)=0\Leftrightarrow\sum_{k=1}^\infty\sum_{i=0}^{k_0-1}
\di^k(\tilde{F}_i)=0
\end{equation}

\begin{figure}[h!]
\begin{center}\label{oog}
\includegraphics[width=6cm]{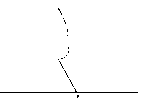}\label{first reduction graph}
\caption{The only graph in $\mathbf{Q}^\infty_{1,1}$.}
\end{center}
\end{figure}
Fix a label $L$. By Lemma \ref{weight diff lemm} one can group together the operators in the differential $\di$ with respect to the number of type I vertices in the respective graphs. Equation (\ref{prohomo}) is then equivalent to a system of equations recovered in the following way: The only operator of weight $\om(L(e_\infty))+\om(\pi)$ is given by the graph $\Gamma$ of Figure 2. Thus the homogeneous equation of weight equal to $k_0+\om(L(e_\infty))+\om(\pi)$ (the highest possible)  in  (\ref{prohomo}) is $\di^1(\tilde{F}_0)=0$. This way we recover the first equation of the system (\ref{system2}). The homogeneous equation in (\ref{prohomo}) of the second highest weight, $k_0+\om(L(e_\infty))+\om(\pi)-1$, is $\di^1(\tilde{F}_1)=0$. For weight $k_0+\om(L(e_\infty))+\om(\pi)-2$, one has $\di^1(\tilde{F}_2)+\di^2(\tilde{F}_0)=0$, which, since all weights $\omega_\Gamma$ for  $\Gamma\in\mathcal{B}^{2i}\cup\mathcal{BW}^{2i}$ are zero, gives $\di^1(\tilde{F}_2)=0$. With the same argument, the equation of weight $k_0+\om(L(e_\infty))+\om(\pi)-3$ is $\di^1(\tilde{F}_3)=0$ and for weight $k_0+\om(L(e_\infty))+\om(\pi)-4$ we get $\di^1(\tilde{F}_4)+\di^3(\tilde{F}_0)=0$, thus recovering the second equation of (\ref{system2}). Similarly, for decreasing weight, one gets $\di^1(\tilde{F}_5)+\di^3(\tilde{F}_1)=~0$, $\di^1(\tilde{F}_6)+\di^3(\tilde{F}_2)=0$ and $\di^1(\tilde{F}_7)+\di^3(\tilde{F}_3)=0$. Inductively, we regroup the terms of $F$ as $F=\sum_iF_{2i}$, where $F_{2i}=\sum_{j=0}^3\tilde{F}_{4i+j}$. \end{proof}

At this point we start discussing a particular choice of weight homogeneous (Lie-) Poisson structure on $\g^\ast$. Let $\g$ be a $k$-dimensional semisimple Lie algebra. Pick an $\mathfrak{sl}_2$-triple $\{e,h,f\}$ and let $n_i$ be the weights of the adjoint action $[h,X_i]=n_iX_i$ on basis elements $X_i\in\g$. Let $S(\g)=\oplus_nS^n(\g)$ be the standard polynomial grading on $S(\g)$. The $\mathrm{ad}h-$ action extends to a derivation on $S(\g)$  and we use the notation $\text{wt}(P)$ for the weight of an element $P\in S(\g)$ with respect to this action. Define $S^n(\g)(i):=\{x\in S^n(\g)/ \mathrm{ad}h(x)=ix\}, i\in\mathbb{Z}$. The Kazhdan grading on $S(\g)$ is $S(\g)=\bigoplus_{n\in\mathbb{Z}}S(\g)[n]$ where $S(\g)[n]:=\bigoplus_{i+2j=n}S^j(\g)(i)$. The Kazhdan degree of a homogeneous element $P\in S(\g)$ will be denoted by $\deg_K(P)$. Let $U_0(\g)\subset U_1(\g)\subset\cdots \subset U_n(\g)\subset\cdots$ be the PBW filtration of $U(\g)$ and set $U_n(\g)(i):=\{x\in U_n(\g)/\mathrm{ad}h(x)=ix\}$. Then the Kazhdan filtration on $U(\g)$ is a $\mathbb{Z}-$ indexed filtration with $F_pU(\g):=\langle x\in U_j(\g)(i)/\;i+2j\leq p\rangle$ (the subspace spanned by all such $x$).
\newtheorem{casap}[con]{Remark}
\begin{casap}\label{kazhdan weight cor}
Set $\om_\g:=(n_1+2,\ldots,n_k+2)$. Then the Kazhdan degree of $P\in S(\g)$ is the weight homogeneous degree $\om_\g(P)$. A direct application of Lemma \ref{kazhdan degree product} is that if $\Gamma\in\mathbf{Q_{1,2}}$, then $\om(B_\Gamma) =-2$. We thus recover as a special case, the fact that the transverse Poisson structure on a slice, has weight homogeneous degree equal to $-2$ (see \cite{DSV}). We should mention that this weight choice is non-trivial; fixing $\om=(1,\ldots,1)$, a linear Poisson structure is weight homogeneous with $\om(\pi)=-1$. It is obvious that for $\om=\om_\g$, Lemma (\ref{k}) simplifies 
to $\sum_{i=1}^n \mathrm{wt}(L(e^1_i))=0$.
\end{casap}

\newtheorem{cas}[con]{Remark}
\begin{cas}\label{degree in H}
Let $\g,\m,\chi,\su$ as in the $W-$ algebra construction. From the proof of Proposition \ref{homo}, we deduce that every element $F\in H^0(\m^\bot,\chi,\su)$ can be written as a finite sum of elements of $S(\su)$ of the form 
\[F=F_0+F_2+F_4+\cdots\]
with $\om_\g(F_{2i})=\om_\g(F_0)-4i$.
\end{cas}
\noindent For the rest of the paper, we fix $\om=\om_\g$ and restrict our grading-related arguments to $\deg_K$.

\noindent\textbf{3.3. A new $W-$ algebra model.}  
\noindent Let $\g_{\epsilon}$ be the Lie algebra over the ring $\mathbb{K}[\epsilon]$ with Lie bracket defined as $[X,Y]_{\epsilon}:=\epsilon[X,Y]$, for $X,Y\in \g$. Set $U_{(\epsilon)}(\g)=U(\g_{\epsilon})$ over $\mathbb{K}[\epsilon]$ and consider the ideal $U_{(\epsilon)}(\g)\m_{\chi+\rho}$ with the notation of Section 2.1.1, where $\rho(H)=\frac{1}{2}\mathrm{Tr}(\mathrm{ad}H)$, $H\in\m$. Denote as $H^0(\g^\ast,\epsilon)$ the $\epsilon-$ reduction algebra corresponding to the data $\su=\g$. Its differential is identically $0$ and it is isomorphic to $S(\g)[\epsilon]$ as a vector space, so $H^0(\g^\ast,\epsilon) \simeq U_{(\epsilon)}(\g)$. Similarly, denote as $H^0(\m^\bot,\g^\ast,\chi,\su,\epsilon)$ the $\epsilon-$ reduction space corresponding to the coisotropic submanifold $\m^\bot$ that is isomorphic to $S(\su)[\epsilon]$ as a vector space and whose differential $\di$ is identically $0$. There is a bimodule structure on $H^0(\m^\bot,\g^\ast,\chi,\su,\epsilon)$: The left $H^0(\g^\ast,\epsilon)-$ module structure is denoted as 
\[\ast_1: H^0(\g^\ast,\epsilon)\times H^0(\m^\bot,\g^\ast,\chi,\su,\epsilon) \longrightarrow H^0(\m^\bot,\g^\ast,\chi,\su,\epsilon)\]
\[ f \ast_1 \rho =f\cdot\rho + \sum_{k=1}^{\infty} \frac{\epsilon ^{k}}{k!}\sum_{\Gamma \in \mathbf{Q}_{k,2}^{(4)}}\overline{\omega}_{\Gamma}B_{\Gamma}(f,\rho)\]
while the right $H^0(\m^\bot,\chi,\su,\epsilon)-$ structure is denoted as  
\[\ast_2: H^0(\m^\bot,\g^\ast,\chi,\su,\epsilon) \times H^0(\m^\bot,\chi,\su,\epsilon) \longrightarrow H^0(\m^\bot,\g^\ast,\chi,\su,\epsilon)\]
\[\rho \ast_2 g =\rho\cdot g + \sum_{k=1}^{\infty} \frac{\epsilon ^{k}}{k!}\sum_{\Gamma\in \mathbf{Q}_{k,2}^{(4)}}\overline{\omega}_{\Gamma}B_{\Gamma}(\rho,g).\]
For more details we invite the reader to check \cite{CF2} and \cite{CT} $\mathcal{x}$ 1.6.
Using these two module structures one defines two differential operators 
\[\overline{T}_1: S(\su)[\epsilon]\longrightarrow  H^0(\m^\bot,\g^\ast,\chi,\su,\epsilon)\]
\[\overline{T}_1(F)=F\ast_1 1\]
\[T_2: H^0(\m^\bot,\chi,\su,\epsilon)\longrightarrow  H^0(\m^\bot,\g^\ast,\chi,\su,\epsilon)\]
\[T_2(G)=1\ast_2 G\]
where the bar over $T_1$ simply denotes the restriction of the operator $\bullet\ast_1 1$ to $S(\su)[\epsilon]\subset H^0(\g^\ast,\epsilon)$. Both $\overline{T}_1$ and $T_2$ are coming from $\mathcal{W}-$ type graphs and so they have constant coefficients. In \cite{BAT},\cite{BAT3} it is proved that there is a non-canonical associative algebra isomorphism,
\begin{equation}\label{theorem}
\overline{\beta}_{\su,(\epsilon)}\circ\partial_{q_{(\epsilon)}^{\frac{1}{2}}}\circ \overline{T}_1^{-1}T_2:\; H^0(\m^\bot,\chi,\su,\epsilon)\stackrel{\sim}{\longrightarrow} \left(U_{(\epsilon)}(\g)/U_{(\epsilon)}(\g)\m_{\chi+\rho}\right)^{\m},
\end{equation}
 where  
\[\overline{\beta}_{\su,(\epsilon)}:\;S(\su)[\epsilon]\longrightarrow U_{(\epsilon)}(\g)/U_{(\epsilon)}(\g)\m_{\chi+\rho}\] 
 is the quotient symmetrization map and 
\[q(Y) := \det_{\g} \left(\frac{\sinh\frac{\mathrm{ad} Y}{2}}{\frac{\mathrm{ad}Y}{2}}\right),\]
for $Y\in\g$. 

 From now on, let $G$ be a connected semisimple Lie group and $\g$ its Lie algebra. Fix a nilpotent element $e\in\g$ and pick $h,f\in\g$ forming an $\mathfrak{sl}_2$- triple with $e$. There exists a $\g-$ invariant bilinear form $(\cdot,\cdot)$ such that $(e,f)=1$. Let $\chi\in\g^\ast$ be defined for all $x\in\g$ by $\chi(x)=(e,x)$. Set  $\g(i):=\{\xi\in\g|[h,\xi]=i\xi\}$ to be the eigenspaces of the $\mathrm{ad}h-$ action. Consider the skew-symmetric form $\omega_{\chi}$ on $\g$ defined by $\omega_{\chi}(\xi,\eta)=\chi([\xi,\eta])$. The restriction $\omega_{\chi}|_{\g(-1)}$ is non-degenerate so one can pick a lagrangian subspace $\mathfrak{l}\subset\g(-1)$. Set $\mathfrak{m}:=\mathfrak{l}\bigoplus\oplus_{i\leq-2}\g(i)$ (so $\chi$ is a character of $\mathfrak{m}$). Set $\mathfrak{m}_{\chi}$ to be the space generated by the elements $\{\xi-\chi(\xi),\;\xi\in\mathfrak{m}\}$ and let $\su$ be such that $\g=\m\oplus\su$.

We prove now that in the $W-$ algebra setup, one can drop the deformation parameter $\epsilon$ and the character $\rho$ in (\ref{theorem}). With the next result we provide a new model of the $W-$ algebra associated to the data $(\g,e)$.

\newtheorem{thea}[con]{Theorem}
\begin{thea}\label{kazhdan theorem}
Let $\g$ be a semisimple Lie algebra, and $\{e,h,f\}$ an $\mathfrak{sl}_2-$ triple. Let $\chi,\m_{\chi}$ be as in the $W-$ algebra setup. There is an associative algebra isomorphism

\begin{equation}\label{main isomorphism}
\overline{\beta}_{\su}\circ\partial_{q^{\frac{1}{2}}}\circ \overline{T}_1^{-1}T_2:\; H^0(\m^\bot,\chi,\su)\stackrel{\simeq}{\longrightarrow} \left(U(\g)/U(\g)\m_{\chi}\right)^{\m}.
\end{equation}
\end{thea}
\begin{proof} The direction $H^0(\m^\bot,\chi,\su)\hookrightarrow\left(U(\g)/U(\g)\m_\chi\right)^{\m}$ works as in the proof of (\ref{theorem}) in \cite{BAT3}, see $\mathcal{x}$ 1.3 of our Introduction. We only note that in the course of this direction, it is proved that if $F\in H^0(\m^\bot,\chi,\su)$, then
\begin{equation}\label{stokes1}
[(m-\chi(m))\ast_11]\ast_2 F=0.
\end{equation}
This is immediate by Lemma 3.3 in \cite{BAT}, 
\begin{equation}\label{h}
(m-\chi(m))\ast_11=0.\end{equation}
The inverse direction uses essentially the same argument with the proof of (\ref{theorem}) in \cite{BAT3}, but we use $\om_\g$ instead of $\deg_\epsilon$. We omit details that can be found in \cite{BAT3}.  
Let $G\in(S(\g)/S(\g)\ast_K\m_{\chi})^{\m}$ be on the vertical axis of the biquantization diagram of $\g^\ast$ and $\m_\chi^\bot$ and let $F=T^{-1}_2\overline{T}_1(G)$ be at the horizontal axis of the diagram. Then if $m-\chi(m),\;m\in\m$, is on the vertical axis, (see Figure (\ref{D0B0})), one has $(m-\chi(m))\ast_{K} G - G\ast_{K} (m-\chi(m)) \in S(\g)\ast_{K} \m_{\chi}$. By a bimodule compatibility relation and (\ref{h})  one gets that $(m-\chi(m))\ast_1(1\ast_2 F)=0$. Thus together with (\ref{stokes1}), we have a Stokes equation 
\begin{equation}\label{stokes}
\sum _{\Gamma}\int^{\infty}_{0}\mathrm{d}\omega_{\Gamma}(s)B_{\Gamma}(F)=0
\end{equation} letting $F$  move as a point $s$ on the horizontal axis. The possible concentrations are:

\noindent\textbf{Interior graphs.} Let $k$ type I vertices and one type II vertex collapse on the horizontal axis.  Dimensional reasons force a possible graph in this concentration to be either of $\mathcal{B}^k-$ or $\mathcal{BW}^k-$ type. Denote as $\alpha$ the edge leaving the concentration and an interior graph as $\Gamma^{\al}_{int}$.

\noindent\textbf{Exterior graphs.} The first possibility is to have a $\mathcal{B}$-type graph receiving at its root the edge $\al$. Then its own $e_{\infty}$ edge derives the function $m-\chi(m)$. Additionally, there can be a superposition of infinitely  many $\mathcal{W}-$ type graphs deriving the concentration. The second possibility is to have only a finite number of $\mathcal{W}-$ graphs deriving the concentration. In this case $\al$ derives $m-\chi(m)$.

\begin{figure}[h!]
\begin{center}\label{D0B0}
\includegraphics[width=5cm]{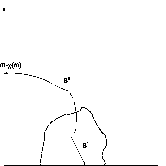}
\caption{\footnotesize The graph corresponding to $B^0_{\Gamma^{\alpha}_{ext}}(B^1_{\Gamma^{\alpha}_{int}})$.}\label{D_0(B_0)}
\end{center}
\end{figure}
\noindent The simplest concentration is the one in Figure \ref{D0B0}, where the interior graph is the one of $\di^{1}(\m^\bot,\chi,\su)$. Let $e_1^1$ be the edge deriving $F$, $e_1^2=\al$ be the edge deriving $m-\chi(m)$ and $L$ a label function. Then by Lemma \ref{weight diff lemm}

\[\deg_K(B^0_{\Gamma^{\alpha}_{ext}}(B^1_{\Gamma^{\alpha}_{int}})(F))=\deg_K(F)-2+\deg_K(L(\alpha)).\]


\noindent In general, suppose that the concentration contains a $\Gamma_{ext}^\al\in\mathcal{B}^m$, $k$ exterior wheel type graphs  $\Gamma_{ext}^{i,\al}\in\mathcal{W}^{r_i}$ and $\Gamma_{int}^\al\in\mathcal{B}^t\cup\mathcal{BW}^t$. Then $\al$ derives the root of $\Gamma_{ext}^\al\in\mathcal{B}^m$. The total Kazhdan degree in this diagram is
\begin{equation}\label{kazhdan theorem 2}
\deg_K(B_{\Gamma_{ext}}(B_{\Gamma_{int}}(F))=
\deg_K(F)+\deg_K(L(\al))-2(t+m)-2\sum_{i=1}^kr_i
\end{equation}
as a straightforward computation shows. The Stokes equation (\ref{stokes}) is  equivalent to
\begin{equation}\label{general}
\sum_{\alpha}\left(\sum_{\Gamma^{\alpha}_{int},\Gamma^{\alpha}_{ext}} \left(B_{\Gamma^{\alpha}_{ext}}(B_{\Gamma^{\alpha}_{int}}(F))\right)\right)=0 \Leftrightarrow \sum_{\alpha}\left(\sum_{\Gamma^{\alpha}_{int},\Gamma^{\alpha}_{ext}}\left(\sum_{l+k+m=0}^{\infty}\left(B^m_{\Gamma^{\alpha}_{ext}}(B^k_{\Gamma^{\alpha}_{int}}(F_l))\right)\right)\right)=0
\end{equation}
for $l=1,\ldots n, \;k=0,\ldots \infty, \;m= 0, \ldots \infty$. Grouping the equations with respect to their weight $\deg_K$,  (\ref{general}) can be rewritten as
\[\sum_{\alpha}B_{\Gamma^{\alpha}_{ext}}(B_{\Gamma^{\alpha}_{int}}(F))=0 \Leftrightarrow \sum_{\alpha}\sum_{l,i,j}B^j_{\Gamma^{\alpha}_{ext}}(B^i_{\Gamma^{\alpha}_{int}}(F_l))=0 \Leftrightarrow \sum_i\di^i(F)=0\Leftrightarrow \di(F)=0.\] \end{proof}

\section{Generators of the reduction algebra.}

\noindent\textbf{4.1. $H^0(\m^\bot,\chi,\su)$ as a deformation of $(S(\g)/S(\g)\m_{\chi})^{\m}$.} 

Let $\g=\q\oplus\m$ and $\chi,\m_{\chi}$ be as in the $W-$ algebra construction. In this section we first prove that  the reduction algebra $H^0(\m^\bot,\chi,\su)$ is a deformation of the space of invariants $(S(\g)/S(\g)\m_{\chi})^{\m}$. We also compute explicitly the generators of $H^0(\m^\bot,\chi,\su)$. For this, we use a specific basis already used by Premet. Let $x_1,\ldots,x_r$ be a basis of $\g_e$ and extend it to a homogeneous basis $x_1,\ldots,x_r,x_{r+1},\ldots,x_m$ of $\p=\sum_{i\geq 0}\g(i)$, with $x_i\in\g(n_i)$ for all $i\in\{1,2,\ldots,m\}$. By \cite{PRE1}, Section 3.1, $[e,\p]=\sum_{i\geq 2}\g(i)$, $\dim \g_e=\dim\g(0)+\dim\g(1)=r$ and  there exist homogeneous elements  $y_{r+1},\ldots,y_m\in\sum_{i\leq -2}\g(i)$, with $y_i\in\g(-n_i-2)$ for all $i\in\{r+1,\ldots,m\}$  such that
\begin{equation}\label{eq1}
\chi([x_i,y_j])=\delta_{ij}\quad(r+1\leq i,j\leq m).\end{equation} In addition, we fix a Witt basis $z_1,\ldots,z_s,z_{s+1},\ldots,z_{2s}$ in $\g(-1)$ such that $[z_{i+s},z_j]=\delta_{ij}f$ and  $[z_i,z_j]=[z_{i+s},z_{j+s}]=0$, for $1\leq i,j\leq s$. Then $\{x_1,\ldots,x_m,z_{s+1},\ldots,z_{2s}\}$ is a basis of $\su$. Let $(\boldsymbol{\alpha,\beta})=(\alpha_1,\ldots,\alpha_m,\beta_1,\ldots,\beta_s)\in\mathbb{N}_0^m\times\mathbb{N}_0^s$ be a multi-index and $x^{\boldsymbol{\alpha}}z^{\boldsymbol{\beta}}=x_1^{\alpha_1}\cdots x_m^{\alpha_m}z_{s+1}^{\beta_1}\cdots z_{2s}^{\beta_s}$ the corresponding element in $S(\su)$. Finally, for $\alpha_i\neq 0$, $\beta_j\neq 0$, we also consider the multi-index $(\boldsymbol{\alpha-e_i,\beta-e_j})=(\alpha_1,\ldots,\alpha_i-1,\ldots,\alpha_m,\beta_1,\ldots,\beta_j-1,\ldots,\beta_s)\in\mathbb{N}^m_0\times\mathbb{N}^s_0$. The $\ast-$ product used throughout this section is the $\ast_{CF}-$ product so we simply denote it as $\ast$.

\newtheorem{lem}[con]{Lemma}
\begin{lem}\label{lem1} Let $P$ be a nonzero Kazhdan-degree-homogeneous polynomial in $(S(\g)/S(\g)\m_{\chi})^{\m}$ of degree $k_0$.  The term in $P$ of minimum polynomial degree is an element of $S(\g_e)$. Namely, $P$ is written in the following form $$P=\sum_{\mathbf{p}=(p_1,\ldots,p_r)}\lambda_{\bf{p}} x_1^{p_1}\cdots x_r^{p_r}+\sum_{\stackrel{(\alpha_{r+1},\ldots,\alpha_m,\beta_1,\ldots,\beta_s)\neq {\bf 0}}{\sum p_i<\sum\alpha_i+\sum\beta_i}} \lambda_{\boldsymbol{\alpha,\beta}} x_1^{\alpha_1}\cdots x_r^{\alpha_r}\cdots x_m^{\alpha_m}z_{s+1}^{\beta_1}\cdots z_{2s}^{\beta_s},$$ with  $\lambda_{\bf{p}},\lambda_{\boldsymbol{\alpha,\beta}}\in \mathbb{K}$.
\end{lem}
\begin{proof} Let $P\in (S(\g)/S(\g)\m_{\chi})^{\m}$ be a nonzero Kazhdan-degree-homogeneous polynomial. With respect to the basis of $\su$ it can be written as
 $$P= \sum_{\deg_K(x^{\boldsymbol{\alpha}}z^{\boldsymbol{\beta}})=k_0}\lambda_{\boldsymbol{\alpha},\boldsymbol{\beta}}x^{\boldsymbol{\alpha}}z^{\boldsymbol{\beta}}.$$ 
Let $(\boldsymbol{\alpha},\boldsymbol{\beta})=(\alpha_1, \ldots,\alpha_r,\alpha_{r+1},\ldots,\alpha_m,\beta_1,\ldots,\beta_s)\in\mathbb{N}^m_0\times\mathbb{N}^s_0$ be an $(m+s)$-tuple such that the polynomial degree of the term $x^{\boldsymbol{\alpha}}z^{\boldsymbol{\beta}}$ is minimum and $(\alpha_{r+1},\ldots,\alpha_m,\beta_1,\ldots,\beta_s)\neq {\bf 0}\in\mathbb{N}^{m-r}_0\times\mathbb{N}^s_0$. If $\alpha_i\neq 0$ for some $i\in\{r+1,\ldots,m\}$, consider the element $y_i\in\g(-n_i-2)\subset\sum_{j\leq -2}\g(j)\subset \m$  defined in (\ref{eq1}). Then we have
\begin{eqnarray*} &  &
[x_1^{\alpha_1}\cdots x_j^{\alpha_j}\cdots x_m^{\alpha_m}z_{s+1}^{\beta_1}\cdots z_{2s}^{\beta_s},y_i]\\ & = & \sum_{j=1}^rx_1^{\alpha_1}\cdots[x_j^{\alpha_j},y_i]\cdots x_r^{\alpha_r}\cdots x_m^{\alpha_m}z_{s+1}^{\beta_1}\cdots z_{2s}^{\beta_s}\\ & + & \sum_{j=r+1}^{m}
\alpha_j\chi([x_j,y_i])x_1^{\alpha_1}\cdots x_r^{\alpha_r}\cdots x_j^{\alpha_j-1}\cdots x_m^{\alpha_m}z_{s+1}^{\beta_1}\cdots z_{2s}^{\beta_s}\\ & = & \sum_{j=1}^rx_1^{\alpha_1}\cdots\alpha_j[x_j,y_i]x_j^{\alpha_j-1}\cdots x_r^{\alpha_r}\cdots x_m^{\alpha_m}z_{s+1}^{\beta_1}\cdots z_{2s}^{\beta_s}\\ & + & \alpha_i x_1^{\alpha_1}\cdots x_r^{\alpha_r}\cdots x_i^{\alpha_i-1}\cdots x_m^{\alpha_m}z_{s+1}^{\beta_1}\cdots z_{2s}^{\beta_s}.
\end{eqnarray*}
For $1\leq j\leq r$, suppose $[x_j,y_i]\in\g(-2)$. Then restricting at $\chi+\mathfrak{m}^\bot$ one has $$[x_j,y_i]=\chi([x_j,y_i])=(e,[x_j,y_i])=-([x_j,e],y_i)=0,$$ since $x_j\in\g_e$. If $[x_j,y_i]\notin\g(-2)$, then $[x_j,y_i]=0$ or $[x_j,y_i]\in\mathfrak{q}$. The second case implies that the term $\sum_{j=1}^rx_1^{\alpha_1}\cdots\alpha_j[x_j,y_i]x_j^{\alpha_j-1}\cdots x_r^{\alpha_r}\cdots x_m^{\alpha_m}z_{s+1}^{\beta_1}\cdots z_{2s}^{\beta_s}$ has greater polynomial degree than the term $\alpha_i x_1^{\alpha_1}\cdots x_r^{\alpha_r}\cdots x_i^{\alpha_i-1}\cdots x_m^{\alpha_m}z_{s+1}^{\beta_1}\cdots z_{2s}^{\beta_s}$. Thus, the nonzero term in $[x_1^{\alpha_1}\cdots x_j^{\alpha_j}\cdots x_m^{\alpha_m}z_{s+1}^{\beta_1}\cdots z_{2s}^{\beta_s},y_i]$ with minimum polynomial degree is $$\alpha_i x_1^{\alpha_1}\cdots x_r^{\alpha_r}\cdots x_i^{\alpha_i-1}\cdots x_m^{\alpha_m}z_{s+1}^{\beta_1}\cdots z_{2s}^{\beta_s}.$$ The $\m$-invariance implies that $\alpha_i=0$, which is a contradiction.

 If $\beta_i\neq 0$ for some $i\in\{1,\ldots,s\}$, we consider the element $z_i\in\g(-1)\subset \m$. A similar argument using $z_i$ instead of $y_i$ results again in contradiction. \end{proof}
\newtheorem{lem2}[con]{Lemma}
\begin{lem2}\label{lem2} Let $Q$ be a nonzero Kazhdan-degree-homogeneous polynomial in $S(\g_e)$ of degree $k_0$. There is a unique polynomial $P\in (S(\g)/S(\g)\m_{\chi})^{\m}$ of Kazhdan degree $k_0$ whose term of minimum polynomial degree is $Q$. 
\end{lem2}
\begin{proof} By Lemma \ref{lem1} every polynomial $P$ in $(S(\g)/S(\g)\m_{\chi})^{\m}$ is written in the form $$P=\sum_{\mathbf{p}=(p_1,\ldots,p_r)}\lambda_{\bf{p}} x_1^{p_1}\cdots x_r^{p_r}+\sum_{\stackrel{(\alpha_{r+1},\ldots,\alpha_m,\beta_1,\ldots,\beta_s)\neq {\bf 0}}{\sum p_i<\sum\alpha_i+\sum\beta_i}} \lambda_{\boldsymbol{\alpha,\beta}} x^{\boldsymbol{\alpha}}z^{\boldsymbol{\beta}},$$ with $\lambda_{\bf{p}},\lambda_{\boldsymbol{\alpha,\beta}}\in \mathbb{K}$. We show that the coefficients $\lambda_{\boldsymbol{\alpha,\beta}}$ are uniquely determined by $\lambda_{\bf{p}}$. 

If $\alpha_i\neq 0$ for some $i\in\{r+1,\ldots,m\}$, consider the element $y_i\in\g(-n_i-2)\subset\sum_{j\leq -2}\g(j)\subset \m$. One then has
\begin{eqnarray*} &  & [\lambda_{\bf{p}} x_1^{p_1}\cdots x_r^{p_r},y_i]=\sum_{j=1}^r p_j\lambda_{\bf{p}} x_1^{p_1}\cdots x_j^{p_j-1}\cdots x_r^{p_r}[x_j,y_i]
\end{eqnarray*} and
\begin{eqnarray*} &  &
[\lambda_{\boldsymbol{\alpha,\beta}}x^{\boldsymbol{\alpha}}z^{\boldsymbol{\beta}},y_i]= \sum_{j=1}^r \alpha_j\lambda_{\boldsymbol{\alpha,\beta}}x^{\boldsymbol{\alpha-e_j}}z^{\boldsymbol{\beta}}[x_j,y_i]+ \alpha_i\lambda_{\boldsymbol{\alpha,\beta}}x^{\boldsymbol{\alpha-e_i}}z^{\boldsymbol{\beta}}.
\end{eqnarray*} We have $\deg_{\q}(x^{\boldsymbol{\alpha-e_j}}z^{\boldsymbol{\beta}}[x_j,y_i])>\deg_{\q}(x^{\boldsymbol{\alpha-e_i}}z^{\boldsymbol{\beta}})$. If $\deg_{\q}(\lambda_{\bf{p}} x_1^{p_1}\cdots x_j^{p_j-1}\cdots x_r^{p_r}[x_j,y_i])=\deg_{\q}(x^{\boldsymbol{\alpha-e_i}}z^{\boldsymbol{\beta}})$, by $\m$-invariance, we deduce that \begin{eqnarray*}
\sum_{j=1}^r p_j\lambda_{\bf{p}} x_1^{p_1}\cdots x_j^{p_j-1}\cdots x_r^{p_r}[x_j,y_i]=-\alpha_i\lambda_{\boldsymbol{\alpha,\beta}}x^{\boldsymbol{\alpha-e_i}}z^{\boldsymbol{\beta}}.
\end{eqnarray*} Consequently, if $\alpha_i\neq 0$, the coefficients $\lambda_{\boldsymbol{\alpha,\beta}}$ are uniquely determined by $\lambda_{\bf{p}}$. Similarly, if $\beta_i\neq 0$ for some $i\in\{1,\ldots,s\}$, it is 
\begin{eqnarray*}
\sum_{j=1}^r p_j\lambda_{\bf{p}} x_1^{p_1}\cdots x_j^{p_j-1}\cdots x_r^{p_r}[x_j,z_i]=-\beta_i\lambda_{\boldsymbol{\alpha,\beta}}x^{\boldsymbol{\alpha}}z^{\boldsymbol{\beta-e_i}},
\end{eqnarray*} so $\lambda_{\boldsymbol{\alpha,\beta}}$ are uniquely determined by $\lambda_{\bf{p}}$. \end{proof}

This lemma implies that there is an algebra isomorphism  $\nu:\,S(\g_e)\longrightarrow (S(\g)/S(\g)\m_{\chi})^{\m}$. Thus, for a basis $x_1,\ldots, x_r$ of $\g_e$, the polynomials $P_i=\nu(x_i)$, $i=1,\ldots,r$  generate $(S(\g)/S(\g)\m_{\chi})^{\m}$. Recall by Remark \ref{degree in H}, that every element in $H^0(\m^\bot,\chi,\su)$ is written as $F=\sum_iF_{2i}$ with $\deg_K(F_{2i})=\deg_K(F_0)-4i$. We show now that $F$ is uniquely determined by the leading term $F_0\in(S(\g)/S(\g)\m_{\chi})^\m$. 
\newtheorem{lem3}[con]{Proposition}
\begin{lem3}\label{lem3}
Let $\displaystyle F=\sum_{i=0}^{[\frac{k_0}{4}]}F_{2i}\in H^0(\m^\bot,\chi,\su)$ with $F_0\in(S(\g)/S(\g)\m_{\chi})^\m$, $\deg_K(F_0)=k_0$. The affine symbol map $\sigma:H^0(\m^\bot,\chi,\su)\longrightarrow (S(\g)/S(\g)\m_{\chi})^\m$ defined by $\sigma(F)=F_0$, is surjective.
\end{lem3}
\begin{proof} 
Let $J=\span\{x_{r+1},\ldots ,x_m,z_{s+1},\ldots,z_{2s}\}$. We have $S(\q)=S(\g_e)\oplus S(\q)J$, where $S(\su)J$ is the ideal of $S(\su)$ generated by $J$. Let $F_2=F'_2+F''_2$ with $F'_2\in S(\g_e)$ and $F''_2\in S(\q)J$. If $F'_2$ is not zero, by Lemma \ref{lem2}, there is an element $\bar{F}'_2$ in $(S(\g)/S(\g)\m_{\chi})^\m$ whose term of minimum polynomial degree is $F'_2$ and the other terms belong to $S(\q)J$. Thus, we may write $F_2=\bar{F}'_2+\bar{F}''_2$, with $\bar{F}'_2\in (S(\g)/S(\g)\m_{\chi})^\m$ and $\bar{F}''_2\in S(\q)J$. Since the invariant term of $F$ is $F_0$, we conclude that $\bar{F}'_2=0$ and $F_2\in S(\q)J$. Thus, $F_2$ is written in the form $$F_2=\sum_{\stackrel{(\alpha_{r+1},\ldots,\alpha_m,\beta_1,\ldots,\beta_s)\neq {\bf 0}}{\deg_K(x^{\boldsymbol{\alpha}}z^{\boldsymbol{\beta}})=k_0-4}} \lambda_{\boldsymbol{\alpha,\beta}}^{k_0-4} x^{\boldsymbol{\alpha}}z^{\boldsymbol{\beta}}.$$ 
By equations (\ref{system}) in Section 3.1 we have that $\di^1(F_2)=-\di^3(F_0)$. Recall by the same section that the graphs appearing in $\di^i$ are either Bernoulli or Bernoulli attached to a wheel and have an edge leaving to infinity that we call $e_\infty$. Let $y_{r+1},\ldots,y_m$ be the elements in $\m$ defined in (\ref{eq1}) and fix $L(e_\infty)=y_j$ for some $j\in\{r+1,\ldots, m\}$. By Lemma (\ref{weight diff lemm}), if the polynomial $\di^3(y_j)(F_0)$ is non-zero, it has Kazhdan degree $k_0-4+\text{wt}(y_j)$. Recall that $\di^1(y_j)(F)=\sum_{i=1}^m[y_j,x_i]\partial_{x_i}(F)$ and that essentially $\di^1(y_j)(F)$ is the bracket $[F,y_j]$, see Remark \ref{inv}. Let now $\alpha_j\neq 0$ for some $j\in\{r+1,\ldots,m\}$. In Kazhdan degree $k_0-4+\text{wt}(y_j)$ we have the following equality of polynomials: 
\begin{eqnarray*} \di^1(\lambda_{\boldsymbol{\alpha,\beta}}^{k_0-4} x^{\boldsymbol{\alpha}}z^{\boldsymbol{\beta}})(y_j) & = & [\lambda_{\boldsymbol{\alpha,\beta}}^{k_0-4} x_1^{\alpha_1}\cdots x_r^{\alpha_r}x_{r+1}^{\alpha_{r+1}}\cdots x_j^{\alpha_j}\cdots x_m^{\alpha_m}z_{s+1}^{\beta_1}\cdots z_{2s}^{\beta_{2s}},y_j]\\ & = & \sum_{k=1}^r\lambda_{\boldsymbol{\alpha,\beta}}^{k_0-4}\alpha_kx_1^{\alpha_1}\cdots x_k^{\alpha_k-1}[x_k,y_j]\cdots x_r^{\alpha_r}x_{r+1}^{\alpha_{r+1}}\cdots x_j^{\alpha_j}\cdots x_m^{\alpha_m}z_{s+1}^{\beta_1}\cdots z_{2s}^{\beta_{2s}}\\ & + & 
\lambda_{\boldsymbol{\alpha,\beta}}^{k_0-4}\alpha_j x_1^{\alpha_1}\cdots x_r^{\alpha_r}x_{r+1}^{\alpha_{r+1}}\cdots x_j^{\alpha_j-1}\cdots x_m^{\alpha_m}z_{s+1}^{\beta_1}\cdots z_{2s}^{\beta_{2s}}.
\end{eqnarray*} 
It is obvious that $\deg_{\q}(x_1^{\alpha_1}\cdots x_k^{\alpha_k-1}[x_k,y_j]\cdots x_m^{\alpha_m}z_{s+1}^{\beta_1}\cdots z_{2s}^{\beta_{2s}})>\deg_{\q}(x_1^{\alpha_1}\cdots  x_j^{\alpha_j-1}\cdots x_m^{\alpha_m}z_{s+1}^{\beta_1}\cdots z_{2s}^{\beta_{2s}})$ for any $k=1,2,\ldots,r$. Thus we equate the last term at the right hand side with the term of $\di^3(y_j)(F_0)$ of equal polynomial degree and compute the coefficients $\lambda_{\boldsymbol{\alpha,\beta}}^{k_0-4}$. If $\beta_j\neq 0$ we use the same arguments for $z_j$,  $j\in\{1,\ldots,s\}$, to compute $\lambda_{\boldsymbol{\alpha,\beta}}^{k_0-4}$. Inductively, using the equation $\di^1(F_{2i})+\di^3(F_{2i-2})+\cdots+\di^{2i+1}(F_0)=0$, we compute the term $F_{2i}$ from $F_0,\ldots,F_{2i-2}$. Hence, for any $F_0\in (S(\g)/S(\g)\m_{\chi})^\m$ there is a unique $F\in H^0(\m^\bot,\chi,\su)$ such that $\sigma(F)=F_0$.
\end{proof}

\noindent In view of Proposition \ref{lem3}, and fixing an invariant  $P\in (S(\g)/S(\g)\m_{\chi})^\m$, we denote as $\tilde{P}$ the element of $H^0(\m^\bot,\chi,\su)$ whose leading term, is $\tilde{P}_0=P$. Obviously, if $k_0\leq 3$, $\tilde{P}=P$.  Before the next Theorem, which is analogous to Theorem 3.4 of \cite{PRE1}, let $B_k$ stand for the sum of bidifferential operators in the $\ast_{CF}-$ product coming from  graphs $\Gamma\in\mathbf{Q}^{(2)}_{k,2}$, taken together with their corresponding coefficient $\omega_\Gamma$. In the proof we ignore possible coefficients in the linear combinations in order to lighten the notation.

\newtheorem{corali}[con]{Theorem}
\begin{corali}\label{paragetai}
Let $\g,\m,\chi$ be as in $\mathcal{x}$ 1.2, $x_1,\ldots,x_r$ a basis of $\g_e$ and $P_i=\nu(x_i)$ the corresponding elements of $(S(\g)/S(\g)\m_{\chi})^\m$. Then $H^0(\m^\bot,\chi,\su)$ is generated by $\tilde{P}_1,\ldots, \tilde{P}_r$.
\end{corali}

\begin{proof} 

By the results of this section, $(S(\g)/S(\g)\m_{\chi})^\m$ is polynomial and isomorphic to $S(\g_e)$. To better illustrate the proof we assume that every element of $(S(\g)/S(\g)\m_{\chi})^\m$ that we consider is written as the product of two generators $P_i,P_j\in (S(\g)/S(\g)\m_{\chi})^\m$ for $i,j\in  \{1,\ldots,r\}$. This is not a loss of generality; if the invariant element is the product of more generators, one can use our procedure and the associativity of the $\ast-$ product to extend it to the other terms as well. Sums of products are handled by linearity of all the operations and operators involved.

 Let $F=\sum_{i=0}^{[\frac{k_0}{4}]}F_{2i}$ with $\deg_K(F_0)=k_0$ be an element of $H^0(\m^\bot,\chi,\su)$. Suppose $F_0=P^0_iP^0_j$ with $P^0_i,P^0_j\in (S(\g)/S(\g)\m_{\chi})^\m$ generating invariants, and $\deg_K(P^0_i)=k_0', \deg_K(P_j^0)=k_0^{''}$. Set 
\[\tilde{P}_i=P^0_i+\sum_{s=1}^{[\frac{k_0^{'}}{4}]}P_i^{2s},\;\;\;\;\tilde{P}_j=P^0_j+\sum_{s=1}^{[\frac{k_0^{''}}{4}]}P_j^{2s}\]
to be the elements of $H^0(\m^\bot,\chi,\su)$ corresponding to $P^0_i,P^0_j$ respectively. Then $\tilde{P}_i\ast\tilde{P}_j=\sum_{m=0}^{[\frac{k_0}{2}]}B_m(\tilde{P}_i,\tilde{P}_j)$ is an element of $H^0(\m^\bot,\chi,\su)$.  
Let $\delta_4$ be the function defined as $\delta_4(k_0)=1$ if $k_0\in 4\mathbb{Z}$ and $\delta_4(k_0)=0$ otherwise.  It is $\deg_K(B_m(\tilde{P}_i,\tilde{P}_j))=k_0-2m$ and $B_m(\tilde{P}_i,\tilde{P}_j)=\sum_{s+k=0}^{[\frac{k_0}{4}]-\delta_4(k_0)}B_m(P_i^{2s},P_j^{2k})$. Then \[\tilde{P}_i\tilde{P}_j=\tilde{P}_i\ast\tilde{P}_j-\sum_{m=1}^{[\frac{k_0}{2}]}B_m(\tilde{P}_i,\tilde{P}_j)\]
and furthermore
\begin{equation}\label{more}
P^0_iP^0_j=\tilde{P}_i\ast\tilde{P}_j-\sum_{s+k=1}^{[\frac{k_0}{4}]-\delta_4(k_0)}P_i^{2s}P_j^{2k}-
\sum_{m=1}^{[\frac{k_0}{2}]}B_m(\tilde{P}_i,\tilde{P}_j)
\end{equation}
with $\deg_K(B_m(P_i^{2s}P_j^{2k}))=k_0-2m-4(s+k)$. From the reduction equations for $F$ and $\tilde{P}_i\ast\tilde{P}_j$, one has
\begin{equation}\label{f}
\di^3(F_0)+\di^1(F_2)=0.
\end{equation}
\begin{equation}\label{two}
\di^3\left(P_i^0P_j^0\right)+\di^1\left(P_i^2P_j^0+P_i^0P_j^2+B_2(P_i^0,P_j^0)\right)=0,
\end{equation}

Subtracting (\ref{two})  from (\ref{f}) we get that $F_2-P_i^2P_j^0-P_i^0P_j^2-B_2(P_i^0,P_j^0)\in (S(\g)/S(\g)\m_\chi)^\m$.
Suppose then that $F_2-P_i^2P_j^0-P_i^0P_j^2-B_2(P_i^0,P_j^0)=P^0_rP^0_t$ where $P^0_t,P^0_r$ are generators of $(S(\g)/S(\g)\m_\chi)^\m$ with degree $\deg_K(P^0_tP^0_r)=k_0-4$. Then as before,
\[ F_2-P_i^2P_j^0-P_i^0P_j^2-B_2(P_i^0,P_j^0)= \tilde{P}_t\ast\tilde{P}_r-\sum_{s+k=1}^{[\frac{k_0}{4}]-\delta_4(k_0)-1}P_t^{2s}P_r^{2k}-\sum_{m=1}^{[\frac{k_0-4}{2}]}B_m(\tilde{P}_t,\tilde{P}_r)\]
with $\deg_K(B_m(P_t^{2s},P_r^{2k}))=k_0-4-2m-4(s+k)$.
Similarly, the second of the reduction equations satisfied by the element corresponding to the invariant $B_1(P_i^0,P_j^0)$ is 
\begin{equation}\label{met}
\di^3(B_1(P_i^0,P_j^0))+\di^1(B_1(P_i^2,P_j^0)+ B_1(P_i^0,P_j^2\})+B_3(P_i^0,P_j^0))=0.
\end{equation}
Since $B_1(P_i^0,P_j^0)\in (S(\g)/S(\g)\m_{\chi})^\m$ and $\deg_K(B_1(P_i^0,P_j^0))=k_0-2$, let $B_1(P_i^0,P_j^0)=P^0_qP^0_\ell$ where $P^0_q,P^0_\ell\in (S(\g)/S(\g)\m_{\chi})^\m$ are generators. As before, one has 
\begin{equation}\label{l}
P^0_qP^0_\ell=\tilde{P}_q\ast\tilde{P}_\ell-\sum_{p+k=1}^{[\frac{k_0-2}{4}]-\delta_4(k_0-2)}P_q^{2p}P_\ell^{2k}-
\sum_{m=1}^{[\frac{k_0}{2}-1]}B_m(\tilde{P}_q,\tilde{P}_\ell)
\end{equation} 
and the reduction equation
\begin{equation}\label{s}
\di^3(P^0_qP^0_\ell)+\di^1(P_q^2P_\ell^0+P_q^0P_\ell^2+B_2(P_q^0,P_l^0))=0.
\end{equation}
Subtracting (\ref{s}) from (\ref{met}), we get that 
\begin{equation}\label{elmt}
B_1(P_i^2,P_j^0)+ B_1(P_i^0,P_j^2)+B_3(P_i^0,P_j^0)-P_q^2P_\ell^0-P_q^0P_\ell^2-B_2(P_q^0,P_\ell^0)
\end{equation}
 is an invariant of Kazhdan degree $k_0-6$. We note here that there is another invariant with that degree, namely $B_1(P_t^0,P_r^0)$. We assume that the element (\ref{elmt}) is written as product of two generators $P^0_m,P^0_n$ of $(S(\g)/S(\g)\m_\chi)^\m$ and continue the process.
Up to this point, 
\[F=\tilde{P}_i\ast\tilde{P}_j+\tilde{P}_t\ast\tilde{P}_r -\sum_{s+k=1}^{[\frac{k_0}{4}]-\delta_4(k_0)-1}P_t^{2s}P_r^{2k} -\sum _{m=1}^{[\frac{k_0-4}{2}]}B_m(\tilde{P}_t,\tilde{P}_r)\]
\[-\sum_{m=1}^{[\frac{k_0}{2}]}\sum_{s+k=0}^{[\frac{k_0}{4}]-\delta_4(k_0)}B_m(P_i^{2s},P_j^{2k})+\sum_{i=2}^{[\frac{k_0}{4}]}F_{2i} -\sum_{s+k=2}^{[\frac{k_0}{4}]-\delta(k_0)}P_i^{2s}P_j^{2k}-B_2(P_i^0,P_j^0)=\]
\[=\tilde{P}_i\ast\tilde{P}_j +\tilde{P}_t\ast\tilde{P}_r + \tilde{P}_q\ast\tilde{P}_\ell+ \tilde{P}_m\ast\tilde{P}_n\;\;+\;\;\text{terms of}\;\; \deg_K\leq k_0-8\]

One can continue to group the remaining terms of $F$ in the expression above by their Kazhdan degree and using the reduction equations of $F, \tilde{P}_i\ast\tilde{P}_j, \tilde{P}_t\ast\tilde{P}_r$ and etc.  These reduction equations will be eventually exhausted. If there are any invariant terms coming from subtraction by parts  with $\deg_K\leq 3$ it will be possible to identify them with some $\tilde{P}=P$. The element $F$ is then written using $\ast-$ products of pairs of elements $\tilde{P}_i$ .
\end{proof}

\end{document}